\documentclass[12pt]{article}
\usepackage{amssymb}

\newcommand{\F}{{\mathbb F}}

\newcommand{\C}{{\mathbb C}}
\newcommand{\EE}{{\mathbb E}}
\newcommand{\SUM}{\raisebox{-0.4ex}{\mbox{\Large $\Sigma$}}}
\newcommand{\eps}{\varepsilon}
\newcommand{\LL}{\Lambda_3}
\newcommand{\PP}{{\mathbb P}}

\newtheorem{theorem}{Theorem}
\newtheorem{lemma}{Lemma}
\newtheorem{corollary}{Corollary}

\title{On the Decay of the Fourier Transform and Three Term Arithmetic
Progressions}

\author{Ernie Croot}

\date{\today}

\begin{document}

\maketitle
\abstract{In this paper we prove a basic theorem of which the following
is a (very weak) corollary:  If 
$f : \F_{p^n} \to [0,1]$ has the property that 
$||\hat f||_{1/3}$ (this is only a quasinorm, because $1/3 < 1$) 
is not too ``large'', and $\EE(f)$ is not too ``small'', 
then there are lots of triples $m,m+d,m+2d$ such that 
$f(m)f(m+d)f(m+2d) > 0$.
If $f$ is the indicator function for some set $S$, then this would be saying that 
the set has many three-term arithmetic progressions.  

In principle this theorem
can be applied to sets having very low density, where $|S|$ is around
$p^{n(1-\gamma)}$ for some small $\gamma > 0$.  

Furthermore, we show that if $g : \F_{p^n} \to [0,1]$ is majorized by $f$, and
$\EE(g)$ is not too ``small'', then in fact
there are lots of progressions $m,m+d,m+2d$ where 
$f(m)g(m+d)f(m+2d) > 0$.}

\section{Introduction}

Suppose that $p$ is a prime number, and $n \geq 1$ is an integer.  Let
$\F$ denote the finite field $\F_{p^n}$ and set $F = |\F|$.  
Suppose that 
$$
f\ :\ \F\ \to\ [0,1].
$$
We will use the expectation operator, defined to be
$$
\EE(f)\ :=\ F^{-1} \SUM_m f(m).
$$
For an $a \in \F$ we will denote the Fourier transform of $f$ at
$a$ as follows
$$\\
\hat f(a)\ =\ \SUM_m f(m) \omega^{a\cdot m},
$$
where $\omega = e^{2\pi i /p}$, and where $a \cdot m$ denotes the
dot product of $a$ and $m$ with respect to the standard $\F_p$ basis
for $\F$.  
\bigskip

Write $\F = \{a_1,...,a_F\}$, where the $a_i$ are ordered so that
$$
|\hat f(a_1)|\ \geq\ |\hat f(a_2)|\ \geq\ \cdots\ \geq\ |\hat f(a_F)|.
$$
For convenience we set $f_i = \hat f(a_i)$; and thus,
$$\\
|f_1|\ \geq\ \cdots\ \geq\ |f_F|.
$$
We also define 
$$
\sigma_i\ =\ \SUM_{i < j \leq F} |f_j|^2,
$$
which is the tail of the spectral $L^2$ norm of $\hat f$.

As a consequence of Parseval we have that if $\EE(f) = \beta$, then
\bigskip

$\bullet$ For all $i=1,...,F$, $\sigma_i \leq \beta F^2$.

$\bullet$ Given $\eps \in (0,1]$, the number of indices $i = 1,...,F$
such that $|f_i| \geq \eps F$ is at most $\beta \eps^{-2}$.
\bigskip

There are certain functions which have a lot fewer ``large'' Fourier 
coefficients as predicted by this second application of Parseval; for example,
suppose that $S$ is a subset of 
$\F$ having $\beta F$ elements, and set 
$$
f(m)\ =\ {1 \over |S|} (S * S)(m).
$$
Then, $f(m)$ is supported on the elements of $S+S$, and clearly
takes on values in $[0,1]$; also, $\EE(f) = \beta$.  Now,
if 
$$
|\hat f(a)|\ =\ {|\hat S(a)|^2 \over |S|}\ \geq\ \eps F,
$$
then 
$$
|\hat S(a)|^2\ \geq\ \eps \beta F^2;
$$
and, by Parseval one can easily show that the number of $a \in \F$
with this property is at most $\eps^{-1}$, which is better than the
$\beta \eps^{-2}$ claimed after the second bullet above (at least for 
fixed $\beta$ and small enough $\eps$).  Furthermore,
the $\sigma_i$ satisfy a sharper inequality than just 
$\sigma_i\ \leq\ \beta F^2$.
In fact, if $i$ is chosen so that $|f_i| \leq \eps F$, then we will
have 
$$
\sigma_i\ \leq\ {|\hat S(a_i)|^2 \over |S|^2} 
\SUM_{i \leq j \leq F} |\hat S(a_j)|^2\ \leq\ \eps F^2. 
$$

If we were to take $f$ to be something like
$$
f(m)\ =\ {1 \over |S|^2} (S * S * S)(m),
$$
we would get even sharper inequalities.
\bigskip

The main theorem of our paper will show that functions like $f$ above 
must always be rich in three-term arithmetic progressions in a certain
sense; actually, it will show even more -- it will show that  there are lots 
of such three-term progressions that pass through dense subsets where
$f$ is positive.

Rather than starting with the statement of this theorem, we will begin by
stating one of its corollaries that is easy to parse.  First, we introduce some
more notation:  Given $f_1,f_2,f_3 : \F \to \C$, define
\begin{eqnarray}
\LL(f_1,f_2,f_3)\ &=&\ \EE_{m,d}(f_1(m)f_2(m+d)f_3(m+2d))\nonumber \\
&=&\ F^{-2} \SUM_{m,d} f_1(m)f_2(m+d)f_3(m+2d).\nonumber
\end{eqnarray}
If all three of our functions $f_1,f_2,f_3$ are the same function $f$, then we
use the abbreviated notation
$$
\LL(f)\ :=\ \LL(f,f,f).
$$
We note that the trivial progressions $m,m,m$ provide the trivial lower bound
$$
\LL(f)\ \geq\ \EE(f)^3 F^{-1}.
$$
We also define the usual norms (and quasinorms for $t \leq 1$)
$$
||f||_t\ =\ \left ( \SUM_a |f(a)|^t \right )^{1/t}.
$$
The corollary alluded to above is as follows.

\begin{corollary} \label{main_theorem_corollary2}  
Suppose $f,g : \F \to [0,1]$, and that 
$$
{\rm For\ all\ } m\in \F,\ f(m)\ \geq\ g(m)\ \geq\ 0;\ {\rm and,\ } 
\EE(f)\ \geq\ \EE(g)\ \geq\ F^{-\theta}.
$$
Then, if 
$$
||\hat f||_{1/3}\ <\ F^{1+\gamma},
$$
we will have that 
$$
\LL(f,g,f),\ \LL(g,f,f)\ \geq\ 10^{-10} p^{-8} F^{-12\theta - 4\gamma}.
$$
\end{corollary}

\noindent {\bf Remark 1.}  It is possible to prove a similar result given constraints
on $||\hat f||_{1/2-\delta}$, for any $\delta > 0$;  however, our method will not 
give good results for quasinorms $1/2$ or higher.
\bigskip

\noindent {\bf Remark 2.}  An example of a function $f$ where this theorem gives
non-trivial results is as follows:  First, let $S$ be a subset of $\F$ having 
$F^{99/100}$ elements.  Then, define
$$
f(m)\ =\ |S|^{-6} (S*S*S*S*S*S*S)(m)    
$$
Note that $f : \F \to [0,1]$, $\EE(f) = \EE(S)$, and $f$ is supported on the 
sumset $S + S + S + S + S + S + S$.  Now,
$$
|\hat f(a)|\ =\ |S|^{-6} |\hat S(a)|^7;
$$
and, using Parseval, we find that the number of places $a$ where 
$$
|\hat S(a)|\ \geq\ 2^{-j} F
$$
is bounded from above by $F^{-1/100} 2^{2j}$.  Thus, 
\begin{eqnarray}
||\hat f||_{1/3}\ &\leq&\ |S|^{-6} \left ( \SUM_{j=0}^\infty 2^{-7j/3} F^{7/3} (F^{-1/100} 2^{2j})
\right )^3
\nonumber \\
&\ll&\ F^{1+3/100}. \nonumber
\end{eqnarray} 
Applying Corollary \ref{main_theorem_corollary2}, it is easy to see that there are
lots of $m$ and $d \neq 0$ such that 
$$
f(m)f(m+d)f(m+2d)\ >\ 0.
$$

Of course, it is fairly easy to prove non-trivial lower bounds for 
$\LL(f)$ without using this corollary
(see \cite{croot}); however, the ideas in the corollary that give these lower 
bounds are different from these other methods (in the Fourier setting the ideas
in \cite{croot} amount to forcing the Fourier transform $\hat f(a)$ to be a positive 
real number at all places $a$; this is quite different from the ideas that lead to
the proof of the above corollary).  
\bigskip

\noindent {\bf Remark 3.}  One way in which this corollary is 
different from others in the theory of arithmetic progressions (e.g. 
\cite{meshulam}), is that it produces 
lower bounds for $\LL(f)$ for when $\EE(f)$ is quite small.  However,
note that the condition that $||\hat f||_{1/3}$ is ``small'' is a very strong 
requirement, only satisfied by certain special ``smooth'' functions,
whereas Meshulam's result \cite{meshulam} holds for arbitrary 
general functions $f$ where $\EE(f) \geq c_p/n$.
\bigskip

The main theorem from which the above corollary follows is:

\begin{theorem}  \label{main_theorem2}
Suppose that $f, g\ :\ \F \to [0,1]$ satisfy
$$
f(m)\ \geq\ g(m)\ \geq\ 0,\ {\rm and\ } \EE(f)\ \geq\ \EE(g)\ \geq\ 
\max(F^{-\theta},\ 8p^{-1/2}k^{-1}).
$$
Further, suppose that 
$$
\sigma_k\ \leq\ \delta^2 F^2.
$$
Then, 
$$
\LL(f,g,f),\ \LL(g,f,f)\ \geq\  p^{-2}{k \choose 2}^{-1} F^{-4\theta}/128 - 9 \delta F^{-2\theta}/8.
$$
\end{theorem}

\noindent {\bf Remark 4.}  One way that one can see how this theorem is
much stronger than the above corollary is as follows:  Say we start with $f$
such that $||\hat f||_{1/3}$ is small enough so that the corollary implies there
are lots of $m,d$ such that $f(m)g(m+d)f(m+2d) > 0$.  Now suppose we change
the value of $f(m)$ at just one place $m$; then, $||\hat f||_{1/3}$ may no longer
be all that small, and the corollary will give only a trivial lower bound
for $\LL(f,g,f)$; however, in a lot of cases,
the change to just one (or in fact many) value of $f(m)$ has little
affect on the value of $\sigma_k$, and so has little affect on the 
conclusion given by the above theorem.
\bigskip

It would be good
if we could have $m,m+d,m+2d$ all three belong to certain special {\it dense} 
subsets of $\F$ (subsets $A$ of $\F$ such that $\EE(f(m)A(m)) > c > 0$); 
however, this appears to be a very difficult and delicate problem to solve, and
would require new ideas in addition to the ones in this paper.

We close the introduction of this paper with the following conjecture, which is
motivated by the above theorem.  We keep the conjecture intentionally vague:
\bigskip

\noindent {\bf Conjecture.}  Suppose $f : \F \to [0,1]$, and $\EE(f) \geq F^{-\theta}$.
Let $k = k(\eps)$ denote the number of places $a \in \F$ where $|\hat f(a)| > \eps F$.
Then, one can obtain a non-trivial bound for $\LL(f)$ purely in terms of 
$\eps$, $k$, and $\theta$.  Basically, what we are asking is a bound of the type
appearing in the above theorem, except that it should not depend on the tail of the
spectral $L^2$ norm of $\hat f$ -- it should only depend on basic information about the
large Fourier coefficients; and, it should give good results when there are only
very few large Fourier coefficients.

\section{Proof of Theorem \ref{main_theorem2} and its corollary}

\subsection{Proof of Corollary \ref{main_theorem_corollary2}}

We first note that from the bound
$$
||\hat f||_{1/3}\ <\ F^{1+\gamma},
$$
we deduce
$$
|f_j|\ <\ {F^{1+\gamma} \over j^3}.
$$
From this it follows that 
$$
\sigma_k\ <\ F^{2 + 2\gamma} \SUM_{j \geq k+1} j^{-6}\ <\ F^{2+ 2\gamma} \int_k^\infty 
x^{-6} dx\ =\ {F^{2+2\gamma} \over 5 k^5}.
$$
Thus, 
$$
\sigma_k\ <\ \delta^2 F^2,\ {\rm for\ } 
\delta\ =\ (2k^{5/2})^{-1} F^\gamma.
$$

From Theorem \ref{main_theorem2} we deduce that 
\begin{eqnarray}
\LL(f,g,f)\ &>&\ p^{-2}{k \choose 2}^{-1} F^{-4\theta}/128 - 9 k^{-5/2} 
F^{-2\theta + \gamma}/16.
\nonumber \\
&\geq&\ p^{-2} k^{-2} F^{-4\theta}/64 - 9k^{-5/2} F^{-2\theta + \gamma}/16. \nonumber
\end{eqnarray}

The value of $k$ which maximizes this last quantity is
$$
k\ =\  2025 p^4 F^{4\theta + 2\gamma},
$$
and it produces the lower bound
$$
\LL(f,g,f)\ \geq\ 10^{-10}p^{-8} F^{-12\theta - 4\gamma}.
$$

\hfill $\blacksquare$

\subsection{Notations and preliminary lemmas}

Let 
$$
A\ =\ \{a_1,...,a_k\}.
$$
denote the set of places corresponding to $f_1,...,f_k$; that is,
$$
f_i\ =\ \hat f(a_i).
$$
Note that because we can have $|f_i| = |f_j|$, the set $A$ is not well
defined; nonetheless,
for the purposes of our proof all we need is that $f_1,...,f_k$ correspond to 
{\it any} set of $k$ largest Fourier coefficients of $f$.
Also, let
$$
B\ :=\ A - A\ =\ \{a-b\ :\ a,b \in A\}.
$$

We seek a subspace $W$ of $\F$ such that 
\bigskip

$\bullet$ At least a quarter of the translates $t \in \F$ (actually, we
only need consider $t \in W^\perp$) satisfy 
\begin{equation} \label{tV}
\SUM_{m \in t + W} g(m)\ \geq\ \EE(g) |W|/2.
\end{equation}

$\bullet$ If $V$ denotes the orthogonal complement of $W$, then
there are no non-zero elements of $B$ that lie in $V$; that is,
\begin{equation} \label{BV0}
B \cap V = \{0\}.
\end{equation}
What this would imply is that all the cosets $a + V$, $a \in A$, are distinct.

$\bullet$ We want $W$ to have small dimension.
\bigskip

We will show that there is a subsapce $W$ satisfying the first two bullets above,
where $|W| = p^{n'}$ (so, $n'$ is the dimension of $W$), where
$$
1 + (\log p)^{-1} \log {k \choose 2}\ \leq\ n'\ <\ 2 + (\log p)^{-1} \log {k \choose 2}.
$$
To this end, we let $S$ denote the set of all subspaces of $\F$ having this 
dimension $n'$.
\bigskip

We begin with a lemma.

\begin{lemma} \label{B_cap_V} 
If we pick a subpace $W \in S$ at random (using uniform measures), 
we will have that if $V = W^\perp$, then 
$$
B\ \cap\ V\ =\ \{0\}.
$$
holds with probability at least $1/2$.
\end{lemma}

\noindent {\bf Proof of the Lemma.} 
Given a random subspace $V$ of codimension $n'$ 
(chosen with the uniform measure),
the probability that some fixed element $b \in B$, $b \neq 0$, lies in 
$V$ will be 
$$
{|V| - 1 \over F - 1} = {p^{n-n'} - 1 \over p^n - 1}.
$$
This follows because $0$ lies in every subspace, and if we eliminate
it, we are left with $F-1$ elements in our field; and, each non-zero
element of the field is just as likely to be in a random subspace $V$
as any other element -- since there are $|V|-1$ non-zero elements of
$V$, this gives the probability $(|V|-1)/(F-1)$.

Thus, since $B$ has at most ${k \choose 2}$ elements, 
the probability that no $b \in B$, $b \neq 0$, lies in $V$ is at least
$$
1\ -\ {k \choose 2} {p^{n-n'} - 1 \over p^n - 1}\ >\ 1\ -\ {k \choose 2}p^{-n'}.
$$
This last quantity exceeds $1/2$ whenever
$$
n'\ \geq\ 1 + (\log p)^{-1} \log {k \choose 2}.
$$
\hfill $\blacksquare$
\bigskip

This Lemma \ref{B_cap_V} is what allows us to produce subspaces $V$ 
satisfying (\ref{BV0}); however, the following lemma will be needed to 
get (\ref{tV}) to hold: 

\begin{lemma} \label{g_lemma}
If
$$
\EE(g)\ >\ 8 p^{-1/2} k^{-1},
$$
then if $t \in \F$ and $W \in S$ are chosen independently at random using
uniform measures, we will have that (\ref{tV}) holds with probability exceeding
$3/4$.
\end{lemma}

\noindent {\bf Proof of the Lemma.}  
The proof of this corollary is via Chebychev's inequality:  Suppose we 
select $t \in \F$ and $W \in S$ independently 
at random using uniform measures.  Define
the random variable
$$
X\ :=\ \SUM_{m \in t + W} g(m).
$$
To prove our corollary it suffices to show that 
$$
{\rm Prob}(| X - |W|\EE(g) |\ >\ |W| \EE(g)/2)\ <\ 1/4.
$$

To prove this using Chebychev, we first consider
$$
\EE(X^2)\ =\ |S|^{-1} F^{-1} 
\SUM_{t \in \F} \SUM_{W \in S} \left ( \SUM_{m \in t + W} g(m) \right )^2.
$$
On expanding out this square, we are left to estimate
$$
\SUM_{m_1, m_2 \in \F} g(m_1)g(m_2) 
\SUM_{ (t,W) \in \F \times S \atop m_1, m_2 \in t + W} 1.
$$
It is easy to see that this equals
$$
\SUM_{m_1, m_2 \in \F} g(m_1) g(m_2) \SUM_{W \in S \atop 
m_1 - m_2 \in W} \SUM_{t \in \F \atop m_1 - t \in W} 1.
$$
(Note that in this final inner sum we get that if $m_1 - t \in W$, then 
$m_2 - t \in W$ as well, because $m_1 - m_2 \in W$.)
Clearly, given $W$ and $m_1$, there are $|W|$ choices for $t \in \F$ 
such that $m_1 - t \in W$; and so, the sum is
$$
\SUM_{m_1, m_2 \in \F} g(m_1) g(m_2) \SUM_{W \in S \atop 
m_1 - m_2 \in W} |W|.
$$

To bound this from above, we consider the case where $m_1 = m_2$
seperate from the case $m_1 \neq m_2$:  The contribution of all $m_1,m_2$
where $m_1 = m_2$ is
$$
\SUM_{m \in \F} g(m)^2 |S| |W|\ \leq\ F |S| |W|\ =\ |S| p^{n + n'}.
$$
The contribution of all unequal pairs $m_1,m_2$ is at most
\begin{eqnarray}
\SUM_{m_1, m_2 \in \F} g(m_1) g(m_2) |S| |W| {|W| -1 \over F - 1}
\ &\leq&\ |S| p^{2n'-n}  \SUM_{m_1,m_2 \in F} g(m_1)g(m_2)  \nonumber \\
&=&\ |S| p^{2n'+n} \EE(g)^2. \nonumber
\end{eqnarray}

So, we deduce that 
\begin{eqnarray}
\EE(X^2)\ &\leq&\ |S|^{-1} F^{-1} \left ( |S| p^{2n'+n} \EE(g)^2 + |S| p^{n'+n} \right )
\nonumber \\
&=&\ p^{2n'} \EE(g)^2 + p^{n'}. \nonumber
\end{eqnarray}

We also have that 
\begin{eqnarray}
\EE(X)\ &=&\ |S|^{-1} F^{-1} \SUM_{W \in S} \SUM_{t \in \F} 
\SUM_{m \in t +W} g(m) \nonumber \\
&=&\ |S|^{-1} F^{-1} \SUM_{m \in \F} g(m) \SUM_{W \in S} \SUM_{t \in \F \atop 
m \in t + W} 1 \nonumber \\
&=&\ |S|^{-1} F^{-1} \SUM_{m \in \F} g(m) \SUM_{W \in S} |W| \nonumber \\
&=&\ p^{n'} \EE(g). \nonumber 
\end{eqnarray}

So, we deduce that 
$$
{\rm Var}(X)\ =\ \EE(X^2) - \EE(X)^2\ \leq\ p^{n'}. 
$$
Chebychev's inequality then gives that 
$$
\PP(|X - \EE(X)| > \EE(X)/2)\ \leq\ {4 {\rm Var}(X) \over p^{2n'} \EE(g)^2}
\ \leq\ {4 \over p^{n'} \EE(g)^2}\ <\ {1 \over 4},
$$
provided
$$
\EE(g)^2\ >\ 64 p^{-1} k^{-2}\ >\ 
16 p^{-1} {k \choose 2}^{-1}\ \geq\ 16 p^{-n'}. 
$$
\hfill $\blacksquare$

A corollary of both Lemmas \ref{B_cap_V} and \ref{g_lemma} is as follows:

\begin{corollary}  \label{main_corollary2} Suppose that 
$$
\EE(g)\ >\ 8 p^{-1/2} k^{-1}.
$$
Then, there exists a subspace $W \in S$ such that

$\bullet$ Equation (\ref{BV0}) holds for $V = W^\perp$; and,

$\bullet$ At least $F/4$ of the translates $t \in \F$ satisfy (\ref{tV}).
\end{corollary}

\noindent {\bf Proof of the Corollary.}  Suppose we select $(t,W) \in \F \times S$
at random using the uniform measure.  Let $E_1$ be the event that 
(\ref{BV0}) holds for $V=W^\perp$, and let $E_2$ be the event that (\ref{tV}) holds.  Then, 
$$
\PP(E_2\ |\ E_1)\ =\ {\PP(E_1,E_2) \over \PP(E_1)}\ \geq\ {\PP(E_1) + \PP(E_2) - 1 
\over \PP(E_1)}.
$$
By Lemma \ref{B_cap_V} we have $\PP(E_1) \geq 1/2$, and by Lemma 
\ref{g_lemma} we have $\PP(E_2) > 3/4$; and so,
$$
\PP(E_2\ |\ E_1)\ >\ 1/4.
$$
It follows that some $W \in S$ has the property that (\ref{BV0}) holds and
that (\ref{tV}) holds for at least $F/4$ translates $t \in \F$.
\hfill $\blacksquare$

\subsection{Construction of the subspace $V$ and the coset $t + W$}

Let $W$ be one of the subspaces described by Corollary \ref{main_corollary2}.
Then, suppose $t \in \F$, and 
define $\alpha := \alpha_t : \F \to \{0,1\}$ to be the indicator function for the
coset $t + W$; that is,
$$\\
\alpha(m)\ =\ \left \{ \begin{array}{rl} 1,\ & {\rm if\ } m \in t + W; \\
                                    0,\ & {\rm if\ } m \not \in t + W.
                  \end{array}\right.
$$
If we let $V = W^\perp$, then the Fourier transform of $\alpha$ is given by
$$\\
\hat \alpha(a)\ =\ \left \{ \begin{array}{rl} |W| \omega^{a \cdot t},\ & {\rm if\ }
a \in V; \\
                                         0,\ & {\rm if\ } a \not \in V.
                       \end{array}\right.
$$

Let 
$$
h(m)\ =\ (f\alpha*V)(m)\ =\ \SUM_{a+b = m} (f\alpha)(a)V(b)\ =\ 
\SUM_{b \in V} (f\alpha)(m-b),
$$
where $V(b)$ denotes the indicator function for $V$.  If $w \not \in W$,
then $\hat h(w) = 0$ (because $\hat V(w) = 0$ in that case); 
however, if $w \in W$, then the Fourier transform
of $h$ is given by
\begin{eqnarray} \label{h_formula2}
\hat h(w)\ &=&\ \widehat{(f\alpha)}(w) \hat V(w)\nonumber \\
&=&\ {1 \over F} (\hat f * \hat \alpha)(w)\hat V(w) \nonumber \\
&=&\ {|V| \over F} \SUM_{u_1 + u_2 = w} \hat f(u_1)\hat \alpha(u_2) \nonumber \\
&=&\ {|V|\cdot |W| \over F} 
\SUM_{u_2 \in V} \hat f(w-u_2) \omega^{u_2 \cdot t}
\nonumber \\
&=&\ \SUM_{v \in V} \hat f(w + v) \omega^{-v \cdot t}. 
\end{eqnarray}
\bigskip

We will now show that there is a choice for $t \in \F$ which guarantees that
the large Fourier spectrum of $h(m)$ is `close' to that of $f(m+t)$:  
First, split $W$ into the union of the sets $W_1$ and $W_2$, where
$W_1$ is the set of all $w \in W$ such that the coset $w + V$ contains
some element of $A$ (which must be unique); $W_2$ is the remaining elements of $W$.  
We use the notations 
$v(x)$ and $w(x)$ to denote
the unique pair of elements of $V$ and $W$, respectively, such that 
$$
x\ =\ v(x) + w(x).
$$
Note that if $x\in A$, then $w(x) \in W_1$.
\bigskip

We seek $t \in \F$ such that the following three things all hold:
\bigskip

$\bullet$  We have that
\begin{equation} \label{first_hf2}
\SUM_{a \in A} 
|\hat h(w(a))\ -\ \omega^{-v(a) \cdot t} \hat f(a)|^2\ \leq\ 4\delta^2 F^2.
\end{equation}

$\bullet$  We have that
\begin{equation} \label{second_hf2}
\SUM_{w \in W_2} |\hat h(w)|^2\ \leq\ 4\delta^2 F^2.
\end{equation}

$\bullet$  Finally, we want that 
\begin{equation} \label{third_bullet}
\SUM_{m \in t + W} g(m)\ \geq\ \EE(g) |W|/2.
\end{equation}
\bigskip

One condition guaranteeing the first two bullets is
\begin{equation} \label{both_sums2}
\SUM_{a \in A} |\hat h(w(a)) - \omega^{-v(a)\cdot t} \hat f(a)|^2
\ +\ \SUM_{w \in W_2} |\hat h(w)|^2\ \leq\ 4\delta^2 F^2.
\end{equation}
From our formula (\ref{h_formula2}) we get that if we sum the 
first sum in (\ref{both_sums2}) over $t \in V$, we get
\begin{eqnarray} \label{two_sums4}
&& \SUM_{t \in V} 
\SUM_{a \in A} |\hat h(w(a)) - \omega^{-v(a) \cdot t} \hat f(a)|^2
\nonumber \\
&&\ \ \ \ \ \ \ =\ 
\SUM_{a \in A} \SUM_{t \in V} \left | \SUM_{v \in V \atop v \neq v(a)}
\hat f(v + w(a)) \omega^{-v \cdot t} \right |^2 \nonumber \\
&&\ \ \ \ \ \ \ =\ 
\SUM_{a \in A} \SUM_{v_1,v_2 \in V \atop v_1,v_2 \neq v(a)} 
\SUM_{t \in V} \hat f(v_1 + w(a)) \overline{\hat f(v_2 + w(a))} 
\omega^{-(v_1 - v_2)\cdot t} \nonumber \\
&&\ \ \ \ \ \ \ =\ 
|V| \SUM_{a \in A} \SUM_{v \in V \atop v \neq v(a)} |\hat f(v+w(a))|^2. 
\end{eqnarray}
If we sum the second sum in (\ref{both_sums2}) over $t \in V$, we get
\begin{eqnarray} \label{two_sums5}
&& \SUM_{w \in W_2} \SUM_{v_1,v_2 \in V} \SUM_{t \in V} \hat f(v_1 + w)
\overline{\hat f(v_2 + w)} \omega^{-(v_1 - v_2)\cdot t} \nonumber \\
&&\ \ \ \ \ \ \ =\ |V| \SUM_{w \in W_2} \SUM_{v \in V} |\hat f(v + w)|^2. 
\end{eqnarray}

The quantities in (\ref{two_sums4}) and (\ref{two_sums5}) sum to 
\begin{equation} \label{V_compare}
|V| \SUM_{a \in \F \atop a \not \in A} |\hat f(a)|^2\ =\ |V| \sigma_k.
\end{equation}
Since the left-hand-side of (\ref{both_sums2}) is invariant under translating $t$ 
by any element of $W$, we deduce that if we extend the sum of the 
left-hand-side of (\ref{both_sums2}) from all $t \in V$ to all $t \in \F$, this sum is
bounded from above by $F \sigma_k$ (instead of $|V|\sigma_k$
as in (\ref{V_compare}) ).  Therefore, 
if we let $T$ denote the set of $t \in \F$ for which (\ref{tV}) holds, then we note that
$|T| \geq F/4$.  We also have that the sum over 
$t \in T \subseteq \F$ of the left-hand-side of (\ref{both_sums2})
is bounded from above by the sum over all $t \in F$, which is 
$F \sigma_k$.  
It follows by simple averaging that there exists $t \in T$ such that 
$$
\SUM_{a \in A} |\hat h(w(a)) - \omega^{-v(a)\cdot t} \hat f(a)|^2
\ +\ \SUM_{w \in W_2} |\hat h(w)|^2\ \leq\ {F \sigma_k \over T}\ \leq\ 
4 \sigma_k;
$$
and so, for this $t \in T$ we will have that both (\ref{first_hf2}) and 
(\ref{second_hf2}) hold; and, trivially, (\ref{third_bullet}) holds by virtue
of the fact that $t \in T$.
\bigskip

\subsection{An $m \in t + W$, $g(m) \geq 0$ is a midpoint of many arithmetic
progressions}

Now select $m \in t + W$ such that 
$$
g(m)\ \geq\ \EE(g)/2.
$$
(By (\ref{tV}) it is obvious such $m$ exists.)
We will show that 
$$
\SUM_d f(m-d) g(m) f(m+d)\ \ {\rm is\ large}.
$$
To do this we just need to show that 
$$
\SUM_d f(m-d) f(m+d)\ \ {\rm is\ large}.
$$
Expressing this in terms of Fourier transforms, we find that it equals
\begin{equation} \label{ff_count}
F^{-1} \SUM_a \hat f(a)^2 \omega^{-2a\cdot m}\ =\ F^{-1} \SUM_{a \in A} 
\hat f(a)^2 \omega^{-2a\cdot m}\ +\ E,
\end{equation}
where the error $E$ satisfies
$$
|E|\ \leq\ F^{-1} \SUM_{j=k+1}^F f_j^2\ =\ F^{-1} \sigma_k\ \leq\ \delta^2 F.
$$
We now compare the final sum in (\ref{ff_count}) with the following:
\begin{equation} \label{compare_sum}
F^{-1} \SUM_{a \in A} \hat h(w(a))^2 \omega^{2v(a)\cdot t - 2a\cdot m}.
\end{equation}

From the Cauchy-Schwarz inequality, we find that these
two sums (\ref{compare_sum}) and the final sum in (\ref{ff_count}) 
differ by at most 
$$
F^{-1} \left ( \SUM_{a \in A} |\hat h(w(a)) - \hat f(a) \omega^{-v(a)\cdot t}|^2 \right )^{1/2}
\left ( \SUM_{a \in A} |\hat h(w(a)) + \hat f(a) \omega^{-v(a)\cdot t}|^2 \right )^{1/2}.
$$
Using (\ref{first_hf2}) and Parseval we find that this is at most 
$$
F^{-1} (2\delta F)(2 F)\ =\ 4 \delta F.
$$
\bigskip

Next, observe that since $m \in t + W$, we have that 
$$
a \cdot m\ =\ v(a)\cdot m + w(a) \cdot m\ =\ v(a) \cdot t + w(a) \cdot m;
$$
and so, the sum in (\ref{compare_sum}) equals
$$
F^{-1} \SUM_{a \in A} \hat h(w(a))^2 \omega^{-2w(a)\cdot m}.
$$

We wish to extend this to a sum over all the elements of $\F$, and to do this we
use the error estimate (\ref{second_hf2}) to deduce that this sum equals
\begin{equation} \label{complete_sum}
F^{-1} \SUM_b \hat h(b)^2 \omega^{-2b\cdot m}\ +\ E',
\end{equation}
where
$$
|E'|\ \leq\ F^{-1} \SUM_{w \in W_2} |\hat h(w)|^2\ \leq\ 4\delta^2 F.
$$
\bigskip

Now, we can interpret the sum in (\ref{complete_sum}) purely in terms of combinatorial
properties of $h$:  The sum equals
$$
\SUM_d h(m-d)h(m+d).
$$
Using the fact that 
$h$ is translation-invariant by elements of $V$, we find that the sum is at least
\begin{eqnarray} 
\SUM_{d \in V} h(m-d)h(m+d)\ =\ h(m)^2 \SUM_{d \in V} 1\ &\geq&\ h(m)^2 |V| \nonumber \\
&=&\ f(m)^2 |V|. \nonumber
\end{eqnarray}
This last equality holds since $h$ and $f$ are equal on the coset $t+W$.
\bigskip

Putting together all our estimates, we find that 
\begin{eqnarray}
\SUM_d f(m-d)f(m+d)\ &\geq&\ f(m)^2 |V|\ -\ 4\delta^2 F\ -\ 4 \delta F\ -\ \delta^2 F \nonumber \\
&\geq&\ g(m)^2 |V|\ -\ 9 \delta F \nonumber \\
&\geq&\ \EE(g)^2|V|/4\ -\ 9 \delta F. \nonumber
\end{eqnarray}

\subsection{From one midpoint to many}

We can repeat the argument in the previous subsection many times for different
values of $m$.  The idea is to reassign $g(m)$ to $0$, to produce the
new function
$$
g_2(x)\ :=\ \left \{ \begin{array}{rl} g(x),\ &\ {\rm if\ } x \neq m; \\
0,\ &\ {\rm if\ } x = m. \end{array} \right.
$$
Then, we find a different $m_2$ where $g_2(m_2) \geq \EE(g_2)/2$, and where
$$
\SUM_d f(m_2-d)f(m_2+d)\ \geq\ \EE(g_2)^2 |V|/4 - 9 \delta F.
$$

Thus, we will produce a sequence of functions $g_1 :=g, g_2,g_3,...,g_r$, and a sequence
of numbers $m_2,m_3,...,m_r$ where 
$r \geq \EE(g)F/2$, each $\EE(g_i) \geq \EE(g)/2$, and $g_i(m_i) \geq \EE(g_i)/2 \geq 
\EE(g)/4$.
We can conclude from this that 
\begin{eqnarray}
\LL(f,g,f)\ &\geq&\ F^{-2} \SUM_{i=1}^r ( \EE(g_i)^2 p^{n-n'}/4 - 9 \delta F) (\EE(g)/4)
\nonumber \\
&\geq&\ F^{-2} (\EE(g)F/2)(\EE(g)/4)(\EE(g)^2 p^{n-n'}/16 - 9 \delta F) \nonumber \\
&=&\ \EE(g)^4 (128 p^{n'})^{-1} - 9\delta \EE(g)^2/8 \nonumber \\
&\geq&\  p^{-2}{k \choose 2}^{-1} F^{-4\theta}/128 - 9 \delta F^{-2\theta}/8. \nonumber
\end{eqnarray}

The proof for $\LL(g,f,f)$ is nearly identical; the only difference is that in this case
we need to bound
\begin{equation} \label{quantity}
F^{-1} \SUM_a \hat f(a) \hat f(-2a) \omega^{a \cdot m},
\end{equation}
instead of 
\begin{equation} \label{quantity2}
F^{-1} \SUM_a \hat f(a)^2 \omega^{-2a\cdot m}.
\end{equation}
The methods we apply above work equally well for (\ref{quantity}) as they do 
for  (\ref{quantity2}).  
\hfill $\blacksquare$


\begin{thebibliography}{999}

\bibitem{croot} E. Croot, I. Ruzsa and T. Schoen, {\it Arithmetic Progressions in 
Sparse Sumsets}, To appear in INTEGERS conference procceedings.

\bibitem{meshulam} R. Meshulam, {\it On subsets of finite abelian groups
with no 3-term arithmetic progressions}, J. Comb. Theory Ser. A. {\bf 71}
(1995), 168-172.


\end{thebibliography}
\end{document}